\documentclass[12pt]{article}
\usepackage{amsmath}
\usepackage{amsfonts}
\usepackage{mathrsfs}
\usepackage{amssymb}

\def\qed{\hfill\mbox{$\Box$}}
\def\Vir{\rm Vir}
\def\a{\alpha}
\def\b{\beta}

\def\d{\delta}

\def\F{\mathbb{F}}
\def\g{\gamma}

\def\l{\lambda}

\def\v{\varphi}

\def\Z{\mathbb{Z}}

\def\ssc{\scriptscriptstyle}

\def\cl{\centerline}

\def\wh{\widehat}
\def\bs{\backslash}
\def\hs{\hspace*}
\def\vs{\vspace*}
\def\ni{\noindent}

\def\BB{{\cal B}}
\def\Derb{{\rm Der\ssc\,}{\cal B}}
\def\adb{{\rm ad\ssc\,}{\cal B}}
\def\Autb{{\rm Aut\ssc\,}{\cal B}}
\def\Intb{{\rm Int\ssc\,}{\cal B}}

\numberwithin{equation}{section}
\newtheorem{theo}{Theorem}[section]

\newtheorem{lemm}[theo]{Lemma}
\newtheorem{prop}[theo]{Proposition}

\usepackage[usenames]{color}

\begin{document}
\cl{\large\bf   Derivations and automorphisms of a Lie algebra of
Block type\footnote{Supported by NSF grant 10825101 of
China\\\indent \ $^*$Corresponding author: chgxia@mail.ustc.edu.cn}}
\vskip5pt

\cl{{Chunguang Xia$^*$, Wei Wang$^\dag$}}

\cl{\small $^{*)}$ \it Department of Mathematics, University of
Science and Technology of China\vs{-4pt}}

\cl{\small\it Hefei 230026, China}

\cl{\small $^{\dag)}$ \it School of Mathematics and Computer
Science, Ningxia University\vs{-4pt}}

\cl{\small\it Yinchuan 750021, China}

\cl{\small\it Email: \ chgxia@mail.ustc.edu.cn, \ \
wwll@mail.ustc.edu.cn}\vs{5pt}

\par\ni {\small{\bf Abstract.} Let $\BB$ be the Lie algebra of Block
type with basis $\{L_{\a,i}|\,\a,i\in\Z, i\geq0\}$ and relations
$[L_{\a,i},L_{\b,j}]=\left((\a-1)(j+1)-(\b-1)(i+1)\right)L_{\a+\b,i+j}.$
In the present paper, the derivation algebra and the automorphism
group of $\BB$ are explicitly described. In particular, it is shown
that the outer derivation space is $1$-dimensional and the inner
automorphism group of $\BB$ is trivial. \vskip5pt

\ni{\bf Key words:} Lie algebras of Block type; derivation;
automorphism.

\ni{\it Mathematics Subject Classification (2000):} 17B05; 17B40;
17B65; 17B68.}

\vskip10pt \ni {\bf 1. \
Introduction}\setcounter{section}{1}\setcounter{theo}{0}

\ni Since a class of infinite dimensional simple Lie algebras was
introduced by Block [B], generalizations of Lie algebras of this
type (usually referred to as {\it Lie algebras of Block type}) have
been studied by many authors (see for example,
[DZ,\,OZ,\,S1,\,S2,\,SZho,\,WT1, WT2,\,X1,\,X2,\,Z,\,ZM]). Thanks to their
relation to the Virasoro algebra, these algebras have attracted more
and more attention in the literature. See for example a survey paper
[S4] on quasifinite representations.

The author in [S2] studied the quasifinite representation of a
family of Lie algebras of this type $\BB(s,G)$ with basis
$\{x^{\a,i}|\,\a\in G,i\in\Z,i\ge0\}$ over an algebraically closed
field $\F$ of characteristic zero and relations
\begin{equation}\label{BsG-block}
[x^{\a,i},x^{\b,j}]=s(\b-\a)x^{\a+\b,i+j}+
((\a-1+s)j-(\b-1+s)i)x^{\a+\b,i+j-1},
\end{equation}
where $G$ is a nonzero additive subgroup of $\F$ and $s=0, 1$.

In [S2], it is pointed out that in case $s=0$, the Lie algebra
$\BB(0,G)$ with $2\in G$ has a nontrivial central extension  induced
by the following 2-cocycle
\begin{equation}\label{BsG-block-2cocycle}
\phi(x^{\a,i},x^{\b,j})=(\a-1)\d_{\a+\b,2}\d_{i,0}\d_{j,0}c,
\end{equation}
where $c$ is a central element. By taking $L_{\a,i}=x^{\a,i+1}$ in
$\BB(0,G)$, we see that the Lie brackets in \eqref{BsG-block} take
the following form
\begin{equation}\label{B-block}
[L_{\a,i},L_{\b,j}]=\left((\a-1)(j+1)-(\b-1)(i+1)\right)L_{\a+\b,i+j}
\mbox{ \ for \ } \a,\b\in G,\,i,j\ge-1.
\end{equation}

In this paper, we focus on  the Lie subalgebra $\BB$ of $\BB(0,\Z)$,
with basis $\{L_{\a,i}|\,\a\in\Z,i\geq0\}$ and the above relations.
The motivation to study this special Block type Lie algebra is
mainly based on a fact that the
central extension, denoted $\wh\BB$, of $\BB$, which is completely
different from $\BB(0,\Z)$ (see (\ref{BsG-block-2cocycle})), is
given by
\begin{equation*}%\label{Wang--}
[L_{\a,i},L_{\b,j}]=\left((\a-1)(j+1)-(\b-1)(i+1)\right)L_{\a+\b,i+j}
+\d_{\a+\b,0}\d_{i,0}\d_{j,0}\frac{\a^3-\a}{6}c,
\end{equation*}
for $\a,\b\in \Z,\,i,j\ge0$, and which contains a subalgebra with
basis $\{L_{\a,0}, c|\,\a\in\Z\}$ isomorphic to the well-known
Virasoro algebra, whereas no central extensions of $\BB(0,\Z)$ can
contain such a subalgebra. Because of this, one may expect that the
representation theory of $\wh\BB$ will be richer and more
interesting than that of $\BB(0,\Z)$ or its central extension.

We realize that the Lie algebra $\BB$ is in fact isomorphic to the
Lie algebra defined in [WT1, WT2] (by regarding $L_{\a-i,i}$ defined
here as $-L_{\a,i}$ defined there). Thus central extensions, modules
of the intermediate series and quasifinite irreducible highest
weight modules of $\BB$ have been considered in [WT1, WT2]. However,
the problem of classification of quasifinite irreducible
$\BB$-modules (which is definitely an important problem in the
representation theory) remains open. It is well understood that the
representation theory of a Lie algebra often depends on its
structure theory. The aim of the present paper is to further study
the structure theory of $\BB$ in order to obtain sufficient
information to give a classification of quasifinite irreducible
$\BB$-modules in the future. In this paper, we first characterize
the structure of the derivation algebra of  $\BB$ and prove that the
outer derivation space or the first cohomology group of $\BB$ with
coefficients in its adjoint module is 1-dimensional (see Theorem
\ref{thm-derivation}). Then we determine the automorphism group of
$\BB$ and show that $\BB$ has no nontrivial inner automorphisms (see
Theorem \ref{thm-automorphiam}).
Finally, we would like to point out that although $\BB$ is
$\Z$-graded with respect to eigenvalues of ${\rm ad}_{L_{0,0}}$, it
is not finitely-generated $\Z$-graded, some classical methods (e.g.,
that in [F]) cannot be applied in our case here.

\vskip15pt \ni {\bf 2. \ Derivations of
$\BB$}\setcounter{section}{2}\setcounter{theo}{0}\setcounter{equation}{0}

\ni Recall that a \emph{derivation} $d$ of the Lie algebra $\BB$ is
a linear transformation on $\BB$ such that
$$
d([x,y])=[d(x),y]+[x,d(y)] \mbox{ \ for \ } x,y\in\BB.
$$
 Denote by $\Derb$ the space of the
derivations of $\BB$ and $\adb$  the space of the \emph{inner
derivations} of $\BB$. It is well known that $\Derb$ forms a Lie
algebra with respect to the commutators of linear transformation of
$\BB$ and $\adb$ is an ideal of $\Derb$. Elements in $\Derb\bs\adb$
are called \emph{outer derivations}. The {\it outer derivation
space} of $\BB$ or the \emph{first cohomology group of $\BB$ with
coefficients in its adjoint module}  is defined by
$$
H^1(\BB)=\Derb /\adb.
$$

Note that $\BB=\oplus_{\a\in\Z}\BB_\a$ is a $\Z$-graded Lie algebra
with $\BB_\a=\mbox{span}\{L_{\a,i}|\,i\in\Z_+\}$. For $\a\in\Z$,
$i\in\Z_+$, we give the following notations
\begin{align*}
&\BB_\a^{[i]}=\mbox{span}\{L_{\a,j}|\,j\leq i\},\ \
 \BB_\a^{(i)}=\mbox{span}\{L_{\a,j}|\,j<i\},\\
&(\Derb)_\a=\{d\in\Derb|\,d(\BB_\b)\subset\BB_{\a+\b} \mbox{\ for\ }
 \b\geq0\}.
\end{align*}
In particular, $\Derb=\oplus_{\a\in\Z}(\Derb)_\a$ is $\Z$-graded.
Obviously, we have a homogeneous derivation of $\BB$ defined by
\begin{equation}\label{equ-der-outer}
d_0:L_{\b,j}\mapsto \b L_{\b,j} \mbox{ \ for \ } \b\in\Z, j\in\Z_+,
\end{equation}
which can be easily verified to be an outer derivation.

\begin{theo}\label{thm-derivation}
The $\Z$-graded derivation algebra
$\Derb=\oplus_{\a\in\Z}(\Derb)_\a$ has the following decomposition:
$$
\Derb=\adb \oplus {\rm D}, \mbox{ \ where \ } {\rm D}={\rm
span}\{d_0\}.
$$
In particular, the first cohomology group of $\BB$ is 1-dimensional,
namely, ${\rm dim\ssc\,}H^1(\BB)=1$.
\end{theo}
\ni{\it Proof.} \ \ Let $d\in\Derb$. The proof of the theorem is
equivalent to proving that $d$ is spanned by ad$_u\in\adb$ for some
$u\in\BB$ and $d_0\in {\rm D}$. This will be done by the following
two lemmas (Lemma \ref{lemma-der-1st} and \ref{lemma-der-2st}).\qed
\vskip5pt

For a fixed integer $\a\in\Z$, consider a nonzero derivation
$d\in(\Derb)_\a$ such that
\begin{equation}\label{equ-der-suppose}
d(\BB^{[j]})\subset\BB^{[i+j]} \mbox{ \ for any \ } j\in\Z_+,
\end{equation}
where $i\in\Z$ is a fixed integer. Using the similar technique as in
[SZho], we assume that the integer $i$  is the minimal one
satisfying (\ref{equ-der-suppose}). Then we can write
\begin{equation}\label{equ-der-supp-modequ}
d(L_{\b,j})\equiv e_{\b,j}L_{\a+\b,i+j}(\mbox{mod }\BB^{(i+j)}),
\end{equation}
where $e_{\b,j}\in\F$ and we adopt the convention that if a notation
is not defined but technically appears in an expression, we always
treat it as zero; for example, $e_{1,0}=0$ if $i<0$ in
(\ref{equ-der-suppose}).

Applying $d$ to
$[L_{\b,j},L_{\g,k}]=\left((\b-1)(k+1)-(\g-1)(j+1)\right)L_{\b+\g,j+k}$,
we have
\begin{equation}\label{equ-der-action}
\begin{split}
&((\a+\b-1)(k+1)-(\g-1)(i+j+1))e_{\b,j}\\
&+((\b-1)(i+k+1)-(\a+\g-1)(j+1))e_{\g,k}\\
&~=((\b-1)(k+1)-(\g-1)(j+1))e_{\b+\g,j+k}.
\end{split}
\end{equation}

\vskip5pt \ni {\bf Claim 1.} \ \ We can assume that $i\in\Z_+$ in
(\ref{equ-der-suppose}).

Otherwise, if $i<0$, then $e_{1,0}=0$ as stated above. Taking $\g=1,
k=0$ in (\ref{equ-der-action}), we have
\begin{equation*}
(\a+\b-1)e_{\b,j} =(\b-1)e_{\b+1,j},
\end{equation*}
which implies that $e_{\b,j}$ does not depend on $j$ for any $\b$.
Letting $j=0$ in (\ref{equ-der-supp-modequ}), we obtain that
$i\geq0$ by the assumption on the minimality of $i$, a
contradiction.

\begin{lemm}\label{lemma-der-1st}
If $\a+i\neq0$ or $\a+i=0$ with $i\neq0$, then $d$ in
$(\ref{equ-der-suppose})$ is an inner derivation.
\end{lemm}

\ni{\it Proof.} \ \ For the case $\a+i\neq0$, taking $\g=k=0$ in
(\ref{equ-der-action}), we have
\begin{equation}\label{equ-der-1-ebj}
(\a+i)e_{\b,j} =((\a-1)(j+1)-(\b-1)(i+1))e_{0,0}.
\end{equation}
Set $u_1=(\a+i)^{-1}e_{0,0}L_{\a,i}\in\BB$ and let $d'=d-\mbox{\rm
ad}_{u_1}$. From (\ref{B-block}) and (\ref{equ-der-1-ebj}) we see
that $d'(L_{\b,j})\in\BB^{(i+j)}$ for  $\b\in\Z$, $j\in\Z_+$. Now by
induction on $i$, one can derive that $d'$ is an inner derivation,
and then $d=d'+\mbox{\rm ad}_{u_1}$ is also an inner derivation.

For the other case $\a+i=0$ with $i\neq0$, we see immediately that
$e_{0,0}=0$ by (\ref{equ-der-1-ebj}).  Applying $d$ to
$[L_{\b-1,j},L_{1,0}]= (\b-2)L_{\b,j}$ and
$\left[L_{\b,j},L_{-1,0}\right] = (\b+2j+1)L_{\b-1,j}$ respectively,
we obtain
\begin{eqnarray}\!\!\!\!\!\!\!\!\!\!\!\!
(\b-i-2)e_{\b-1,j}+((\b-2)(i+1)+i(j+1))e_{1,0}\!\!\!&=&\!\!\!
  (\b-2)e_{\b,j}, \label{equ-der-2-L10}\\\!\!\!\!\!\!\!\!\!\!\!\!
(\b+i+2j+1)e_{\b,j}\!+\!((\b-1)(i+1)+(i+2)(j+1))e_{-1,0}\!\!\!&=&\!\!\!
  (\b+2j+1)e_{\b-1,j}.\label{equ-der-2-L-10}
\end{eqnarray}
In particular, taking $\b=j=0$ in (\ref{equ-der-2-L10}), we see that
\begin{equation}\label{equ-der-2-e-10e10}
e_{-1,0}+e_{1,0}=0.
\end{equation}
Multiplying (\ref{equ-der-2-L10}) by $\b+2j+1$,
(\ref{equ-der-2-L-10}) by $\b-i-2$, and then adding both results
together, we obtain $i(i+2j+3)(e_{\b,j}-(\b+j)e_{1,0})=0$ by
(\ref{equ-der-2-e-10e10}), which implies for $i\ne0$ that
\begin{equation}\label{equ-der-2-ebj}
e_{\b,j}=(\b+j)e_{1,0} \mbox{ \ for \ } \b\in\Z, j\in\Z_+.
\end{equation}
Set $u_2=-\frac{1}{i+1}e_{1,0}L_{-i,i}\in\BB$ and let
$d''=d-\mbox{\rm ad}_{u_2}$. By (\ref{B-block}) and
(\ref{equ-der-2-ebj}), we obtain that $d''(L_{\b,j})\in\BB^{(i+j)}$
for $\b\in\Z$, $j\in\Z_+$. As in the first case, by induction on
$i$, we see that $d''$ is an inner derivation, thus $d$ is also an
inner derivation.\qed

\begin{lemm}\label{lemma-der-2st}
If $\a=i=0$, then $d$ in {\rm(\ref{equ-der-suppose})} can be written
as $d=\mbox{\rm ad}_u+\l d_0$ for some $u\in\BB$ and $\l\in\F$.
\end{lemm}

\ni{\it Proof.} \ \ Now the equations (\ref{equ-der-2-L10}) and
(\ref{equ-der-2-L-10}) can be simplified as
$$
\begin{array}{rclr}
\hs{122pt}(\b-2)(e_{\b-1,j}+e_{1,0}-e_{\b,j})&=&0,&(\ref{equ-der-2-L10}')\\[5pt]
(\b+2j+1)(e_{\b-1,j}-e_{-1,0}-e_{\b,j})&=&0.&\hs{122pt}(\ref{equ-der-2-L-10}')
\end{array}
$$
We claim that
\begin{equation}\label{equ-der-3-ebj''}
e_{\b,j}=\b e_{1,0}+e_{0,j} \mbox{ \ for \ } \b\in\Z, j\in\Z_+.
\end{equation}
In fact, if $\b\neq2$, then $e_{\b,j}=e_{1,0}+e_{\b-1,j}$ by
(\ref{equ-der-2-L10}$'$). By induction on $\b$, one can easily
obtain that
\begin{eqnarray}\label{equ-der-3-ebj'}
 e_{\b,j}=
\left\{
\begin{aligned}
 &\b e_{1,0}+e_{0,j}&& \mbox{if} \ \ \b\leq1,\\
 &(\b-2)e_{1,0}+e_{2,j}&& \mbox{if} \ \  \b\geq3.
\end{aligned} \right.
\end{eqnarray}
If $\b=2$, then
$e_{2,j}=e_{1,j}-e_{-1,0}=e_{1,j}+e_{1,0}=2e_{1,0}+e_{0,j}$ by
(\ref{equ-der-2-L-10}$'$), (\ref{equ-der-2-e-10e10}) and the first
case of (\ref{equ-der-3-ebj'}) respectively. This, together with
(\ref{equ-der-3-ebj'}), gives the claim.

On the other hand, the equation (\ref{equ-der-action}) can be
rewritten as
$((\b-1)(k+1)-(\g-1)(j+1))(e_{\b,j}+e_{\g,k}-e_{\b+\g,j+k})=0$.
Substituting \eqref{equ-der-3-ebj''} in this formula gives
$$
((\b-1)(k+1)-(\g-1)(j+1))(e_{0,j}+e_{0,k}-e_{0,j+k})=0.
$$
Then $e_{0,j+k}=e_{0,j}+e_{0,k}$ by arbitrariness of $\b$ or $\g$.
By induction on $j$, one can derive that $e_{0,j}=je_{0,1}$, which,
together with \eqref{equ-der-3-ebj''}, gives
\begin{equation}\label{equ-der-3-ebj}
e_{\b,j}=\b e_{1,0}+j e_{0,1} \mbox{ \ for \ } \b\in\Z, j\in\Z_+.
\end{equation}
Set
$$
\bar{d}=d+\mbox{\rm ad}_{u_3}-(e_{1,0}-e_{0,1})d_0,
$$
where $u_3=e_{0,1}L_{0,0}\in\BB$ and $d_0$ is defined by
(\ref{equ-der-outer}). Applying $\bar{d}$ to the formula
$[L_{0,0},L_{\b,j}]=-(\b+j)L_{\b,j}$, using (\ref{equ-der-3-ebj}),
we obtain that $\bar{d}(L_{\b,j})\in\BB^{(j)}$ for $\b\in\Z$,
$j\in\Z_+$. By Lemma \ref{lemma-der-1st}, $\bar{d}$ is an inner
derivation, and then $d=\mbox{\rm ad}_u+(e_{1,0}-e_{0,1})d_0$ for
some $u\in\BB$. This completes the proof.\qed

\vskip15pt \ni {\bf 3. \ Automorphisms of
$\BB$}\setcounter{section}{3}\setcounter{theo}{0}\setcounter{equation}{0}
\vskip5pt

\ni  An element $S\in\BB$ is called
\begin{itemize}\parskip0pt
\item[(i)] {\it $\rm ad$-locally finite} if for any
given $v\in \BB$, the subspace ${\rm Span}\{{\rm ad}_S^m\cdot
v|\,m\in\Z_+\}$ of $\BB$ is finite dimensional,
\item[(ii)] {\it $\rm ad$-locally nilpotent} if for any given
$v\!\in\!\BB$, there exists some $N\!>\!0$ such that ${\rm
ad}_S^N\cdot v\!=\!0$.
\end{itemize}
Denote by $\Autb$ the \emph{automorphism group} of $\BB$, and
$\Intb$ the \emph{inner automorphism group} of $\BB$, namely, the
subgroup of $\Autb$, generated by ${\rm exp}^{{\rm ad}_x}$ for $\rm
ad$-locally nilpotent elements $x$'s.

In this section, we first prove that $\BB$ does not have a nonzero
locally nilpotent element, thus the inner automorphism group of
$\BB$ is trivial. Next we construct three kinds of outer
automorphisms of $\BB$, and then completely characterize the
structure of the automorphism group of the Lie algebra $\BB$.

\begin{lemm}\label{lemma-aut-local-finite}
Up to scalars, $L_{0,0}$ is  the unique locally finite element of
$\BB$. Furthermore, $\BB$ does not have a nonzero locally nilpotent
element, thus the inner automorphism group of $\BB$ is trivial.
\end{lemm}

\ni{\it Proof.} \ \ Take any locally finite element
$S=\sum_{(\a,i)\in I_S}\l_{\a,i}L_{\a,i}$ of $\BB$, where $I_S$ is a
finite subset of $\Z\times\Z_+$. First, suppose that there exists
$\l_{\a,i}\neq0$ for some $\a<0$. Take the minimal $\a_0<0$ such
that there exists some $i$ with $\l_{\a_0,i}\neq0$, and then choose
$i=i_0$ to be the maximal one satisfying this condition. By
rescaling $S$, we may suppose
$$
S=L_{\a_0,i_0}+\sum_{\a>\a_0\,\,{\rm or}\,\atop\a=\a_0,
i<i_0}\l_{\a,i}L_{\a,i},
$$
and in this case we say that $S$ has the {\it minimal term}
$L_{\a_0,i_0}$.
 Recall that $[L_{\a_0,i_0},L_{\b,j}]=F^j_\b
L_{\a_0+\b,i_0+j}$, where we use the following notation
$$
F^j_\b:=(\a_0-1)(j+1)-(\b-1)(i_0+1).
$$
If $\a_0+i_0\geq0$ (or $>0)$, we can choose big (or small) enough
$\b_0$ and suitable $j_0$ such that $F^{j_0}_{\b_0}<0$ ( or $>0$)
and
$$
F^{j_0+ki_0}_{\b_0+k\a_0}=F^{j_0}_{\b_0}-k(\a_0+i_0)<0\,\ (\mbox{or
}>0) \mbox{ \ for all \ } k\in\Z_+,
$$
which implies that ${\rm ad}_S^k(L_{\b_0,j_0})$, with minimal terms
$L_{\b_0+k\a_0,j_0+ki_0}$, are linear independent for all $k$, i.e.,
$S$ is not $\rm ad$-locally finite. Hence $\l_{\a,i}=0$ for all
$\a<0$. Similarly, we can also show that $\l_{\a,i}=0$ for all
$\a>0$.

Now we can rewrite $S=\sum_{i\in I'_S}\l_{0,i} L_{0,i}$, where
$I'_S$ is a finite subset of $\Z_+$. If there exists $\l_{0,i}\neq0$
for some $i>0$, then similarly take $i_0>0$ to be the maximal one
and assume that
$$
S=L_{0,i_0}+\sum_{i<i_0}\l_{0,i} L_{0,i}.
$$
Now $[L_{0,i_0},L_{\b,j}]=G^j_\b L_{\b,i_0+j}$, where
$G^j_\b=-(j+1)-(\b-1)(i_0+1)$. One can take big enough $\b_0$ and
some $j_0$ satisfying
$$
G^{j_0+ki_0}_{\b_0}=G^{j_0}_{\b_0}-k i_0<0 \mbox{ \ for \ }
k\in\Z_+,
$$
which also contradicts our assumption. So $\l_{0,i}=0$ for all
$i>0$, and thus $S=\l_{0,0} L_{0,0}$ for some $\l_{0,0}\in\F$,
namely, $L_{0,0}$ is up to scalars the unique locally finite element
of $\BB$.

Note that any locally nilpotent element must be locally finite
element by definition. Since ${\rm
ad}^N_{L_{0,0}}L_{\a,i}=-N(\a+i)L_{\a,i}\ne0$ for any $N>0$ if
$\a+i\ne0$, we know that the locally finite element $L_{0,0}$ is not
locally nilpotent. Hence the above statement implies that $\BB$ does
not have a nonzero locally nilpotent element, and then the inner
automorphism group of $\BB$ is trivial. \qed\vskip5pt

Recall that the centerless Virasoro algebra $\Vir$ with basis
$\{L_\a|\,\a\in\Z\}$ is defined by the commutation relations:
$[L_\a,L_\b]=(\b-\a)L_{\a+\b}$ for $\a,\b\in\Z$. We review a known
result about the structure of the automorphism group of Virasoro
algebra. It can also be regarded as a corollary of Theorem 2.3 in
[SZha]. \vskip5pt

\begin{prop}\label{prop-aut-vir}
\begin{itemize}\parskip0pt
\item[{\rm(1)}] For any $\mu\in\F^*$, the
following map is an automorphism of $\Vir$.
$$
\chi_\mu:{\Vir\rightarrow\Vir},\ \ L_\a\mapsto\mu^\a L_\a \mbox{ \
for any \ } \a\in\Z.
$$
\item[{\rm(2)}] For any $s\in\{\pm1\}\cong\Z/2\Z$, the following
map is an automorphism of $\Vir$.
$$
\chi'_s:{\Vir\rightarrow\Vir},\ \ L_\a\mapsto s L_{s\a} \mbox{ \ for
any \ } \a\in\Z.
$$
\item[{\rm(3)}] ${\rm Aut}(\Vir)\cong\F^*\rtimes\Z/2\Z.$
\end{itemize}
\end{prop}

Motivated by the above, one can define the following three kinds of
maps:
\begin{eqnarray}
&\v_\mu: & \BB\rightarrow\BB \ \ \  L_{\a,i}\mapsto \mu^\a L_{\a,i}; \nonumber\\
&\v'_\nu: & \BB\rightarrow\BB \ \ \  L_{\a,i}\mapsto \nu^i L_{\a,i}; \nonumber\\
&\rho_\xi: & \BB\rightarrow\BB \ \ \  L_{\a,i}\mapsto \xi
 L_{\xi(\a+i)-i,i},\nonumber
\end{eqnarray}
where $\mu,\nu\in\F^*=\F\bs\{0\}$ and $\xi\in\{\pm1\}$. One can
easily check that they are all (outer) automorphisms of $\BB$.
Furthermore, we have the following facts.

{\rm (1)} $\{\v_\mu|\,\mu\in\F^*\}\cong\F^*$ is a subgroup of
$\Autb$, where $\v_{\mu_1}\v_{\mu_2}=\v_{\mu_1\mu_2}$ for
$\mu_1,\mu_2\in\F^*$.

{\rm (2)} $\{\v'_\nu|\,\nu\in\F^*\}\cong\F^*$ is a subgroup of
$\Autb$, where $\v'_{\nu_1}\v'_{\nu_2}=\v'_{\nu_1\nu_2}$ for
$\nu_1,\nu_2\in\F^*$.

{\rm (3)} $\{\rho_\xi|\,\xi=-1,1\}\cong\Z/2\Z$ is a subgroup of
$\Autb$.

\begin{prop}\label{prop-aut-Sub-Vir}
Let ${\cal V}={\rm Span}\{L'_\a\,|\,\a\in\Z\}$ be a subalgebra of
$\BB$, which is isomorphic to the centerless Virasoro algebra, i.e.,
$[L'_a,L'_\b]=(\a-\b)L'_{\a+\b}$. Suppose  $L'_0\in\F L_{0,0}$. Then
$L'_\a\in\F L_{\a,0}$ for all $\a\in\Z$.
\end{prop}

\ni{\it Proof.} \ \ By rescaling $L'_0$, we can suppose
$L'_0=L_{0,0}$. Let $0\ne\a\in\Z$. Write $L'_{\a}=\sum_{(\b,j)\in
J_\a}\mu_{\b,j}L_{\b,j}$, where $J_\a$ is a finite subset of
$\Z\times\Z_+$. Then
$$
-\a\sum_{(\b,j)\in J_\a}\mu_{\b,j}L_{\b,j}=-\a L'_\a=[L'_0,L'_\a]=
\Big[L_{0,0},\sum_{(\b,j)\in J_\a}\mu_{\b,j}L_{\b,j}\Big]
=-\sum_{(\b,j)\in J_\a}(\b+j)\mu_{\b,j}L_{\b,j},
$$
which implies that $\mu_{\b,j}=0$ if $\b+j\ne \a$. Hence we can
rewrite
$$
L'_{\a}=\sum_{j\in J'_\a}\l_{\a,j} L_{\a-j,j}, \mbox{ where
}\l_{\a,j}=\mu_{\a-j,j},\ J'_\a=\{j\,|\,(\a-j,j)\in
J_a\}\subset\Z_+.
$$
Then
\begin{eqnarray}\label{equ-aut-prop}
 2\a L_{0,0}&=&[L'_\a,L'_{-\a}]=\Big[\sum_{i\in J'_\a}\l_{\a,i}
L_{\a-i,i},\sum_{j\in J'_{-\a}}\l_{-\a,j} L_{-\a-j,j}\Big]\nonumber\\
&=&\sum_{(i,j)\in J'_\a\times J'_{-\a}}\!\!\!\!(i+j+2)\a\l_{\a,i}
\l_{-\a,j} L_{-(i+j),i+j}.
\end{eqnarray}
Let $i_0={\rm max}\{i\,|\,i\in J'_{\a},\,\l_{\a,i}\ne0\}$,
$j_{0}={\rm max}\{j\,|\,j\in J'_{-\a},\,\l_{-\a,j}\ne0\}$. If
$i_0+j_0>0$, then the right-hand side of \eqref{equ-aut-prop}
contains the nonzero term $(i_0+j_0+2)\a\l_{\a,i_0} \l_{-\a,j_0}
L_{-(i_0+j_0),i_0+j_0}$, which is not in $\F L_{0,0}$. Thus
$i_0=j_0=0$ (since $i_0,j_0$ are non-negative), in particular
$L'_\a\in\F L_{\a,0}$. \qed\vskip5pt

\begin{lemm}\label{lemma-aut-a0}
Let $\tau\in\Autb$, then $\tau(L_{\a,0})=\xi\mu^\a L_{\xi\a,0}$ for
some $\mu\in\F^*$, and $\xi\in\{\pm1\}$.
\end{lemm}

\ni{\it Proof.} \ \ Suppose $\tau\in\Autb$. Let
$L'_\a=\tau(L_{\a,0})$ for $\a\in\Z$. Since ${\cal N}={\rm
Span}\{L_{\a,0}\,|\,\a\in\Z\}$ is the centerless Virasoro algebra,
we see that ${\cal V}={\rm Span}\{L'_\a\,|\,\a\in\Z\}$ is a
subalgebra isomorphic the centerless Virasoro algebra. Furthermore,
since $L_{0,0}$ is up to scalars the unique ${\rm ad}$-locally
finite element in $\BB$, we must have $L'_0=\tau(L_{0,0})\in\F
L_{0,0}$. So Proposition \ref{prop-aut-Sub-Vir} implies $\tau({\cal
N})={\cal V}={\cal N}$. Now the result follows from Proposition
\ref{prop-aut-vir}.\qed

\begin{lemm}\label{lemma-aut-0i}
Let $\tau\in\Autb$, then $\tau(L_{0,i})=\xi\nu^i L_{(\xi-1)i,i}$ for
some $\nu\in\F^*$, and $\xi\in\{\pm1\}$.
\end{lemm}

\ni{\it Proof.} \ \ Assume
\begin{equation}\label{equ-aut-L0i}
\tau(L_{0,i})=\sum_{(p,q)\in J_i}\nu_{p,q}L_{p,q}~~\mbox{for some}~~
\nu_{p,q}\in\F,
\end{equation}
where  $J_i$ is some finite subset of $\Z\times\Z_+$. Applying
$\tau$ to the equation $[L_{0,0},L_{0,i}]=-i L_{0,i}$, we get
$$
\sum_{(p,q)\in J_i}(i-\xi(p+q))\nu_{p,q}L_{p,q}=0,
$$
which implies that $\nu_{p,q}=0$ if $p\ne\xi i-q$. Then
(\ref{equ-aut-L0i}) can be rewritten as
\begin{equation}\label{equ-aut-L0i'}
\tau(L_{0,i})=\sum_{q\in J'_i}\l_{i,q}L_{\xi i-q,q}, \mbox{ where }
\l_{i,q}=\nu_{\xi i-q,q},\ \ J'_i=\{q\,|\,(\xi i-q,q)\in J_i\}.
\end{equation}
Applying $\tau$ to $[L_{-1,0},[L_{1,0},L_{0,i}]]=-2(i+1)L_{0,i}$,
using Lemma \ref{lemma-aut-a0}, we obtain
$$
\sum_{q\in J'_i}(q-i+1)(q+i+2)\l_{i,q}L_{\xi i-q,q}=
2(i+1)\sum_{q\in J'_i}\l_{i,q}L_{\xi i-q,q},
$$
which then implies that $(q-i)(q+i+3)\l_{i,q}=0$, and thus
$\l_{i,q}=0$ if $q\ne i$. Thus we can rewrite (\ref{equ-aut-L0i'})
as
\begin{equation*}
\tau(L_{0,i})=\xi \nu_i L_{(\xi-1)i,i} \mbox{ \ for some \ }
\nu_i\neq0.
\end{equation*}
Finally, applying $\tau$ to the relation
$[L_{0,i},L_{0,1}]=(i-1)L_{0,i+1}$, we obtain $\nu_{i+1}=\nu_i
\nu_1$, which implies $\nu_i=\nu^i$, where $\nu=\nu_1$, and the
lemma follows.\qed

\begin{theo}\label{thm-automorphiam}
Let $\tau\in\Autb$, then there exist some
$\mu,\nu\in\F^*,\,\xi\in\{\pm1\}$ such that $$\tau(L_{\a,i})=\xi
\mu^\a \nu^i L_{\xi(\a+i)-i,i} \mbox{ \ for \ }
\a,i\in\Z,\,i\in\Z_+.$$ In particular, $
\Autb\cong(\F^*\times\F^*)\rtimes\Z/2\Z. $
\end{theo}

\ni{\it Proof.} \ \ Let $\tau\in\Autb$, by Lemma \ref{lemma-aut-a0}
and \ref{lemma-aut-0i}, we have $\tau(L_{\a,0})=\xi\mu^\a
L_{\xi\a,0}$ and $\tau(L_{0,i})=\xi\nu^i L_{(\xi-1)i,i}$ for some
$\mu,\nu\in\F^*$ and $\xi\in\{\pm1\}$. Applying $\tau$ to the
equation $[L_{\a,0},L_{0,i}]=(\a(i+1)-i)L_{\a,i}$ gives
$$
(\a(i+1)-i)\left(\tau(L_{\a,i})-\xi\mu^\a \nu^i
L_{\xi(\a+i)-i,i}\right)=0.
$$
Thus the result holds if $\a(i+1)\neq i$. Assume $\a(i+1)=i$, which
implies $\a=i=0$ since $i\in\Z_+$. In this case, we have the result
by Lemma \ref{lemma-aut-a0}.\qed

\small
\vs{15pt}\par\ni {\bf References} \vs{7pt} \baselineskip=-2pt
\parskip=-0.03truein
\parindent=2.0cm
\def\hang{\hangindent\parindent}
\def\textindent#1{\indent\llap{[#1]\enspace}\ignorespaces}
\def\re{\par\hang\textindent}
\small\def\bf{}

\re{B} R. Block, On torsion-free abelian groups and Lie algebras,
   \textit{Proc. Amer. Math. Soc.} \textbf{9} (1958) 613--620.

\re{DZ} D. Dokovic, K. Zhao, Derivations, isomorphisms and
  second cohomology of generalized Block algebras,
  \textit{Algebra Colloq.} \textbf{3} (1996) 245--272.

\re{F} R. Farnsteiner, Derivations and central extensions of
finitely generated graded Lie algebras, \textit{J. Algebra}
\textbf{118} (1988) 33--45.

\re{OZ} J.M. Osborn, K. Zhao, Infinite-dimensional Lie algebras
  of generalized Block type,  \textit{Proc. Amer. Math. Soc.}
  \textbf{127} (1999) 1641--1650.

\re{S1} Y. Su, Quasifinite representations of a Lie algebra of Block
  type, \textit{J. Algebra} \textbf{276} (2004) 117--128.

\re{S2} Y. Su, Quasifinite representations of a family of Lie
  algebras of Block type, \textit{J.~Pure Appl. Algebra} \textbf{192}
  (2004) 293--305.

\re{S3} Y. Su,  2-Cocycles on the Lie algebras of generalized
  differential operators, \textit{Comm. Algebra} \textbf{30} (2002)
  763--782.

\re{S4} Y. Su, Quasifinite representations of some Lie algebras
related to the Virasoro algebra, in: Advanced Lectures in Math.,
vol. 8, 2009, pp. 213--238.

\re{SZha} Y. Su, K. Zhao, Generalized Virasoro and super-Virasoro
algebras and modules of the intermediate series, \textit{J. Algebra}
\textbf{252} (2002) 1--19.

 \re{SZho} Y. Su, J. Zhou, Structure of the Lie algebras related to
those of Block, \textit{Comm. Algebra} \textbf{30} (2002)
  3205--3226.

\re{WT1} Q. Wang, S. Tan, Quasifinite modules of a Lie algebra
  related to Block type, \textit{J. Pure Appl. Algebra} \textbf{211}
  (2007) 596--608.

\re{WT2} Q. Wang, S. Tan, Leibniz central extension on a Block Lie
algebra, \textit{Algebra Colloq.} \bf{14} (2007) 713-720.

\re{X1} X. Xu, Generalizations of Block algebras,
  \textit{Manuscripta Math.} \textbf{100} (1999) 489--518.

\re{X2} X. Xu, Quadratic conformal superalgebras,
  \textit{J. Algebra} \textbf{224} (2000) 1--38.

\re{Z} K. Zhao, A class of infinite dimensional simple Lie algebras,
  \textit{J. London Math. Soc.} (2) \textbf{62} (2000) 71--84.

\re{ZM} L. Zhu, D. Meng, Structure of degenerate Block algebras,
  \textit{Algebra Colloq.} \textbf{10} (2003) 53--62.

\end{document}